\documentclass[reqno, 12pt]{amsart}

\usepackage[]{_pascal/_pascal}

\hypersetup{
    pdftitle={Treewidth versus Hadwiger number: Graph classes with dominated balanced separators are quasi-linearly treewidth-Hadwiger bounded},
    pdfauthor={Clément Dallard, J. Pascal Gollin, Andrea Munaro},
}

\newcommand{\had}{\ensuremath{\mathsf{had}}}
\newcommand{\tw}{\ensuremath{\mathsf{tw}}}

\newcommand{\halftoone}{\ensuremath{[\nicefrac{1}{2},1)}}

\newcommand{\Free}[2]{\ensuremath{\mathcal{F}_{\mathsf{#1}}(#2)}}

\begin{document}

\author[Dallard]{Cl\'{e}ment Dallard$^\ast$}
\address[$\ast$]{ Department of Informatics, University of Fribourg, Fribourg, Switzerland}
\email{\tt clement.dallard@unifr.ch}

\author[Gollin]{J.~Pascal Gollin$^\dagger$}
\address[$\dagger$]{FAMNIT, University of Primorska, Koper, Slovenia}
\email{\tt pascal.gollin@famnit.upr.si}

\author[Munaro]{Andrea Munaro$^\ddagger$}
\address[$\ddagger$]{Department of Mathematical, Physical and Computer Sciences, University of Parma, Italy}
\email{\tt andrea.munaro@unipr.it}

\title[Treewidth versus Hadwiger number: Dominated balanced separators]{Treewidth versus Hadwiger number:\\ Graph classes with dominated balanced separators are quasi-linearly treewidth-Hadwiger bounded}

\date{\today}


\subjclass[2020]{}

\begin{abstract}
    In a hereditary graph class~$\mathcal{G}$, the treewidth of a graph in~$\mathcal{G}$ is bounded by a function of the size of a largest clique minor of the graph if and only if the class~$\mathcal{G}$ excludes some planar graph as an induced minor.  
    Campbell et al.~(arXiv 2024)~conjectured that such a bounding function can always be chosen to be linear. 
    We verify this conjecture for several classes of graphs (including even-hole-free graphs) and show that if a hereditary graph class~$\mathcal{G}$ has dominated binary-weighted balanced separators, then the bounding function can be taken to be in~${O(h \log h)}$, where $h$ is the Hadwiger number.
    We say that~$\mathcal{G}$ has dominated binary-weighted balanced separators if, for every~${G \in \mathcal{G}}$ and every~$\mathsf{w} \colon V(G) \to \{0,1\}$, there is a set~${D \subseteq V(G)}$ of bounded size so that deleting the closed neighbourhood~$N_G[D]$ of $D$ leaves each component with at most a fixed fraction of the total weight of~$G$. 
    Our result reduces the quasi-linear weakening of the conjecture of Campbell et al.~to a conjecture of Gartland and Lokshtanov stating that every graph class excluding a planar graph as an induced minor admits dominated binary-weighted balanced separators. 
\end{abstract}

\maketitle

\section{Introduction}

One of the most prominent parameters in structural graph theory is the \textsl{treewidth}\footnote{See \zcref{sec:prelims} for the precise definition of this concept.} $\tw(G)$ of a graph~$G$, which is often described as the standard measure of how similar~$G$ is to a tree.
Treewidth provides an ingredient for the structural description of many graph classes and allows for many algorithmic applications. 
The \emph{Hadwiger number~$\had(G)$} of a graph~$G$ is the largest integer~$h$ such that~$G$ contains~$K_h$, the complete graph on~$h$ vertices, as a minor. 
Although~${\had(G) \leq \tw(G)+1}$, there is no converse bound in general: 
planar grids have arbitrarily large treewidth and Hadwiger number at most four. 
It is therefore natural to ask which hereditary\footnote{A graph class is \emph{hereditary} if it is closed under induced subgraphs.} graph classes admit a bound on treewidth in terms of Hadwiger number, and how quickly such a bound must grow. 
We say a hereditary graph class~$\mathcal{G}$ is \emph{${(\tw,\had)}$-bounded} if there is a function~$f$ (called a \emph{$(\tw,\had)$-bounding function of $\mathcal{G}$}) such that~${\tw(G) \leq f(\had(G))}$ for every graph~${G \in \mathcal{G}}$. 

For a hereditary graph class~$\mathcal{G}$ to be ${(\tw,\had)}$-bounded it is clearly necessary that~$\mathcal{G}$ does not contain planar minor models of arbitrarily large grids. 
This necessitates the exclusion of some planar graph~$H$ as an \textsl{induced minor}. 

We say that a graph~$H$ is an \emph{induced minor} of a graph~$G$ if~$H$ is isomorphic to a graph that can be obtained from~$G$ by deleting vertices and contracting edges (and deleting any resulting parallel edges). 
We say that for a graph~$H$, a graph~$G$ is \emph{$H${\dash}induced-minor-free} if~$G$ does not contain~$H$ as an induced minor, and we say that a hereditary graph class $\mathcal{G}$ \emph{excludes~$H$ as an induced minor} if each~${G \in \mathcal{G}}$ is $H${\dash}induced-minor-free. 
Let~$\Free{im}{H}$ be the class of all graphs excluding~$H$ as an induced minor. 

It is a somewhat folkloric consequence of the Grid Theorem~\cite{RobertsonSeymour:1986:GM5} and Flat Wall Theorem~\cite{RobertsonSeymour:1995:GM13} of Robertson and Seymour that the condition of excluding a planar graph~$H$ as an induced minor is also sufficient for $(\tw,\had)$-boundedness. 
Related variants of such a result have appeared multiple times in recent years \cite{AboulkerAKST:2021:evenholefreetw, AbrishamiACHS:2024:ISTD7, Bonnet:2025:SparseInducedSubgraphs}. 

Campbell et al.~\cite{CampbellDDFGHHWWY:2024} initiated the systematic study of ${(\tw,\had)}$-bounded graph classes and were the first to attempt to optimise the order of the bounding function. 

\begin{theorem}[{Campbell et al.~\cite{CampbellDDFGHHWWY:2024}}]
    \label{thm:tw-had-boundedness}
    A hereditary graph class is $(\tw,\had)$-bounded if and only if it excludes a planar graph as an induced minor.
    For each such graph class there exists a $(\tw,\had)$-bounding function in~$O(h^9(\log h)^{O(1)})$. 
\end{theorem}

Their proof relies on the best known bounds for Robertson and Seymour's Grid Theorem that were obtained by Chuzhoy and Tan~\cite{ChuzhoyTan:2021:PolyGrid}, but they conjecture that a linear bound always suffices. 
We call a graph class \emph{linearly} or \emph{quasi-linearly} $(\tw,\had)$-bounded if there exists a $(\tw,\had)$-bounding function in~${O(h)}$ or~${O(h (\log h)^{O(1)})}$, respectively. 
Note that all quasi-linear bounds we obtain in this work are actually in~${O(h \log h)}$. 

\begin{conjecture}[{Campbell et al.~\cite{CampbellDDFGHHWWY:2024}}]
    \label{conj:main}
    Every ${(\tw,\had)}$-bounded hereditary graph class is linearly ${(\tw,\had)}$-bounded.
\end{conjecture}

This conjecture was verified for several graph classes, 
including circle graphs~\cite{HickingbothamIMW:2024:Circlegraphs}, 
bounded rank perturbations of circle graphs, 
outer-string graphs, 
${(g,c)}$-outer-string graphs\footnote{A graph is a ${(g,c)}$-outer-string graph if it is an intersection graph of curves on a surface~$\Sigma$ obtained from a surface of Euler genus~$g$ by deleting the interior of~$c$ disjoint closed discs, and each curve has an endpoint on the boundary of~$\Sigma$.}, 
and ordered graphs excluding an induced cross~\cite{CampbellDDFGHHWWY:2024}. 

In this work, we extend the number of classes for which this conjecture holds, as well as identify a large number of classes which are quasi-linearly $(\tw,\had)$-bounded. 

To this end, we take a closer look at the necessary condition of excluding a planar graph~$H$ as an induced minor. 
It is well known that grids are universal graphs for planar graphs with respect to both the minor relation (see, e.g.,~\cite[(1.5)]{RobertsonST:1994:quicklyplanar}) and the induced minor relation (see, e.g.,~\cite[Theorem~12]{CampbellDDFGHHWWY:2024}). 
So let us first think about strengthening the condition of excluding large grids as induced minors. 

Consider the analogously defined notion of \emph{$(\tw,\omega)$-boundedness}, where $\omega$ denotes the \textsl{clique-number} of a graph. 
This notion was introduced by Dallard et al.~\cite{DallardMS:2024:TWvsOmega1,DallardMS:2024:TWvsOmega2,DallardMS:2024:TWvsOmega3}. 
Clearly, every $(\tw,\omega)$-bounded graph class is $(\tw,\had)$-bounded. 
In contrast to \zcref{conj:main}, Chudnovsky and Trotignon~\cite{ChudnovskyTrotignon:2024:twvsomegacounterexamples} showed that for any function~$g$ there is a $(\tw,\omega)$-bounded graph class for which every $(\tw,\omega)$-bounding function is in~$\Omega(g)$. 
Dallard et al.~\cite{DallardMS:2024:TWvsOmega2} also introduced the stronger notion of \textsl{tree-independence number}\footnote{See \zcref{sec:prelims} for the precise definition of this concept.}. 
They showed that a graph class with bounded tree-independence number is $(\tw,\omega)$-bounded~\cite{DallardMS:2024:TWvsOmega2}, while the converse fails by a result of Chudnovsky and Trotignon~\cite{ChudnovskyTrotignon:2024:twvsomegacounterexamples}. 

By using a classic result of Duchet and Meyniel~\cite{DM82} (see \zcref{thm:duchetmeyniel}), it can be immediately deduced that classes of bounded tree-independence number are linearly $(\tw,\had)$-bounded. 

\begin{theorem}
    [store=thmatw]
    \label{thm:atw}
    Let~$c$ be a positive integer and let~$\mathcal{G}$ be a hereditary graph class such that every graph in~$\mathcal{G}$ has tree-independence number at most~$c$. 
    Then ${f(h) \coloneqq (2c-1)h-1}$ is a $(\tw,\had)$-bounding function of~$\mathcal{G}$. 
\end{theorem}

Another approach to strengthening the condition of excluding large grids is via the concept of \textsl{balanced separators}. 
If we delete a small set of vertices from a large grid, then there is a component of the remaining graph that contains the majority of the vertices of the grid. 
Dvo\v{r}\'ak and Norin~\cite{DvorakNorin:2019:BalancedSeps} showed that a graph class has bounded treewidth if and only if every subgraph of every graph in the class admits a balanced separator of uniformly bounded size. 
Equivalently, there is a constant~$k$ such that, for every graph~$G$ in the class and every set~${X \subseteq V(G)}$, there is a set ${S \subseteq V(G)}$ of size at most~$k$ such that each component of~${G-S}$ contains at most~${\abs{X}/2}$ vertices of~$X$.
By the Grid Theorem of Robertson and Seymour~\cite{RobertsonSeymour:1986:GM5}, these conditions are also equivalent to excluding some fixed planar graph as a minor. 
Gartland and Lokshtanov 
conjectured a variant of this characterisation with respect to induced minors, where they replace the bounded size of these sets by a bounded number of vertices dominating each set. 

Let us introduce the following notation. 
We write~$\halftoone$ for the half-open interval of all reals~${x}$ with~${\nicefrac{1}{2} \leq x < 1}$. 
For a positive integer~$d$ and~$c \in \halftoone$, a hereditary graph class~$\mathcal{G}$ has \emph{$d$-dominated binary-weighted $c$-balanced separators} if, for every~${G \in \mathcal{G}}$ and every~$\mathsf{w} \colon V(G) \to \{0,1\}$, there is a set~${S \subseteq V(G)}$ of size at most~$d$ such that~${\mathsf{w}(V(C)) \leq c \cdot \mathsf{w}(V(G))}$ for each component~$C$ of~${G - N[S]}$. 
We say~$\mathcal{G}$ has \emph{$d$-dominated $c$-balanced separators} if the above is only required for constant weightings~$\mathsf{w}$. 

It is again not hard to see that deleting the closed neighbourhood of a small set of vertices in a large grid leaves a component of the remaining graph that contains the majority of the vertices of the grid. 
On the other hand, Gartland and Lokshtanov (see~\cite{Gartland:2023:Thesis}) conjectured that for every planar graph~$H$ there exists a positive integer~$d$ such that the class of $H${\dash}induced-minor-free graphs has $d$-dominated $\nicefrac{1}{2}$-balanced separators. 
Their conjecture is equivalent to the following binary-weighted version; we sketch a short proof  of this equivalence based on arguments from~\cite{Gartland:2023:Thesis} in the appendix. 

\begin{conjecture}
    \label{conj:balancedseps}
    For every planar graph~$H$ there exists a positive integer~$d$ and a~${c \in \halftoone}$ such that $\Free{im}{H}$ has $d$-dominated binary-weighted $c$-balanced separators. 
\end{conjecture}

This discussion leads to our main result. 
Replacing the condition of ``excluding a planar graph as an induced minor'' by the (possibly) stronger condition of ``having $d$-dominated binary-weighted $c$-balanced separators'' is sufficient to prove quasi-linear $(\tw,\had)$-boundedness. 

\begin{theorem}
    [store=thmquasi]
    \label{thm:main}
    Let~$d$ be a positive integer, let~${c \in \halftoone}$, and let~$\mathcal{G}$ be a hereditary graph class with $d$-dominated binary-weighted $c$-balanced separators. 
    Then~$\mathcal{G}$ is $(\tw,\had)$-bounded with a $(\tw,\had)$-bounding function~$f(h) \in O(h \log h)$. 
\end{theorem}

 Now, \zcref{thm:tw-had-boundedness,thm:main} imply that the weakening of \zcref{conj:main} to quasi-linear $(\tw,\had)$-boundedness is implied by \zcref{conj:balancedseps}. 

\begin{theorem}
    \label{thm:main+conj}
    \zcref{conj:balancedseps} implies that every $(\tw,\had)$-bounded hereditary graph class admits a $(\tw,\had)$-bounding function~$f(h) \in O(h \log h)$. 
\end{theorem}
 
\zcref{conj:balancedseps} is verified for various graphs~$H$ (see \zcref{thm:dombalancedseps}), including paths, cycles, fans, and wheels. 
For a positive integer~$t$, let~$P_t$ and~$C_t$ denote the path and the cycle with $t$ vertices, respectively, let~$F_t$ denote the \emph{$t$-fan}, that is, the graph obtained from~$P_t$ by adding a universal vertex, and let~$W_t$ denote the \emph{$t$-wheel}, that is, the graph obtained from~$C_t$ by adding a universal vertex. 
In particular we get the following corollary.

\begin{corollary}
    [store=corquasi]
    \label{cor:main+conj}
    For each integer~${t \geq 3}$, 
    the classes~$\Free{im}{P_t}$, $\Free{im}{C_t}$, $\Free{im}{F_t}$, $\Free{im}{W_t}$ are $(\tw,\had)$-bounded with $(\tw,\had)$-bounding functions in~$O(h \log h)$. 
\end{corollary}

\smallskip

As mentioned above, in some cases, our methods yield linear $(\tw,\had)$-boundedness. 

An \emph{even-hole-free graph} is a graph~$G$ not containing any even cycle as an induced subgraph. 
Aboulker et al.~\cite{AboulkerAKST:2021:evenholefreetw} showed that even-hole-free graphs are $(\tw,\had)$-bounded, while they are known not to be~$(\tw,\omega)$-bounded as shown by Sintiari and Trotignon~\cite{SintiariTrotignon:2021:LayeredWheels}. 
We verify \zcref{conj:main} for even-hole-free graphs and, more generally, \textsl{$C_4$-free odd-signable graphs}.\footnote{See \zcref{sec:prelims} for the precise definition of this class.}

\begin{theorem}[store=thmevenholefree]
    \label{thmevenholefree}
    The class of~$C_4$-free odd-signable graphs and the class of even-hole-free graphs are linearly $(\tw,\had)$-bounded. 
\end{theorem}

Our final result is that if in addition to a planar graph we exclude a fixed complete bipartite graph~$K_{t,t}$ as an induced minor, then dominated binary-weighted balanced separators suffice to prove linear $(\tw,\had)$-boundedness. 

\begin{theorem}[store=thmHandKtt]
    \label{thmHandKtt}
    For every positive integer~$t$ and every planar graph~$H$, if there exists a positive integer~$d$ and~${c \in \halftoone}$ such that the class~$\mathcal{G} \coloneqq \Free{im}{H} \cap \Free{im}{K_{t,t}}$ has $d$-dominated binary-weighted $c$-balanced separators, then~$\mathcal{G}$ is linearly $(\tw,\had)$-bounded. 
\end{theorem}

\subsection{Overview of the proofs}

Throughout this overview, $c$ and~$d$ are fixed.
Our proofs have two main ingredients. 
The first is a reduction proven in \zcref{sec:method} (as \zcref{thm:mainmethod}).
A \emph{$(k,2k)$-block} is a set of $2k$ vertices that cannot be separated by deleting fewer than~$k$ vertices.
Adapting an argument of Chudnovsky et al.~\cite{ChudnovskyGHLS:2026:EvenHoleFree}, we use lean tree-decompositions to show that, in the absence of such a block, the separator hypothesis yields balanced separators of size~$O(k)$ for every binary weighting. 
The characterization of treewidth by weighted balanced separators then gives~${\tw(G)=O(k)}$. 

For our second main ingredient, we show that $(k,2k)$-blocks for large~$k$ force large $K_h$-minors in the respective settings, which we will do in \zcref{sec:blocks2cliqueminors}. 
Applying this with ${h=\had(G)+1}$ and combining it with \zcref{thm:mainmethod} proves the respective theorems, which we elaborate on in \zcref{sec:applications}. 

To obtain this for \zcref{thmevenholefree}, we mainly build on the aforementioned theorem of Duchet and Meyniel~\cite{DM82} (see \zcref{thm:duchetmeyniel}). 
Using their result, in the absence of large clique minors, we can also reduce large $(2kh,4kh)$-blocks to \emph{strong $(k,2k)$-blocks}, that is, independent sets of~$2k$ vertices between each pair of which we find~$k$ pairwise internally vertex-disjoint and anticomplete paths. 
Note that these conditions are not required to hold for the paths between different pairs of vertices. 
Using a result of Chudnovsky et al.~\cite{ChudnovskyCFL:2025} guaranteeing that large strong blocks always admit a~$K_{t,t}$ or a fixed planar graph~$H$ as an induced minor, we are able to prove \zcref{thmHandKtt}. 

For \zcref{thm:main}, the second main ingredient is \zcref{thm:balancedblock}, which shows that, under the separator hypothesis, a $(k,2k)$-block forces a $K_h$-minor when~$k$ is a sufficiently large multiple of~${h\log h}$. 
The proof of \zcref{thm:balancedblock} is provided in \zcref{sec:balancedblock} and proceeds as follows. 

Starting from a $(k,2k)$-block~$X$, \zcref{lem:separator-cut-family} uses repeated applications of the separator hypothesis to produce~$\Theta(k)$ pairwise disjoint sets~$D_i$, each of size at most~$d$. 
Their union~$Z$ is small enough that a subset~${M \subseteq X \setminus Z}$ of size $\Theta(k)$ retains sufficient inseparability in~${G-Z}$.
For each~$i$, either~$D_i$ has~$\Omega(k)$ neighbours in~$M$, or there is a bipartition of~$M$ with both sides
of size~$\Omega(k)$ such that every path in~${G-Z}$ joining opposite sides meets the neighbourhood of~$D_i$. 

Next, \zcref{lem:ARV-dense-minor} converts this structure into a minor of average degree~${\Omega(k/\sqrt{\log k})}$.
If the first alternative holds for~$\Omega(k)$ indices, the adjacencies between the sets~$D_i$ and~$M$ already give a sufficiently dense subgraph. 
Otherwise, we encode the vertices of~$M$ by unit vectors recording their sides in the remaining bipartitions. 
The geometric separation theorem of Arora, Rao, and Vazirani (\zcref{thm:ARV-separation}) gives two subsets~${U, V \subseteq M}$, each of size~$\Theta(k)$, such that every pair in~${U \times V}$ is separated by~${\Omega(k/\sqrt{\log k})}$ bipartitions.
Menger's Theorem supplies~$\Theta(k)$ pairwise vertex-disjoint paths from~$U$ to~$V$ in~${G-Z}$, each adjacent to that many distinct sets~$D_i$. 
Contracting these paths and retaining the vertices of~$Z$ as singleton branch sets gives the required minor.

Equipped with a minor of average degree~${\Omega(k/\sqrt{\log k})}$, we use the well-known Kostochka--Thomason theorem (\zcref{thm:kostochka-thomason}) to obtain a large clique minor. 
To force a $K_h$-minor, it suffices that~${k/\sqrt{\log k}}$ be a sufficiently large multiple of~${h\sqrt{\log h}}$, and choosing~$k$ to be a sufficiently large multiple of~${h\log h}$ ensures this, thus concluding the proof of \zcref{thm:balancedblock}.

\section{Preliminaries}
\label{sec:prelims}

All graphs in this paper are finite and simple. 
All logarithms are base~$2$. 
For standard graph theory notation not introduced here we refer the reader to~\cite{Diestel:2025:GraphTheory}. 

Let~$G$ be a graph with vertex set~$V(G)$ and edge set~$E(G)$. 

For a set~${X \subseteq V(G)}$, we write $N_G(X)$ for the \emph{open neighbourhood of~$X$}, that is, the set ${\{ v \in V(G) \colon \exists w \in X, vw \in E(G) \} \setminus X}$, 
and we write $N_G[X]$ for the \emph{closed neighbourhood of~$X$}, that is, the set ${\{ v \in V(G) \colon \exists w \in X, vw \in E(G) \} \cup X}$. 
If~${X = \{v\}}$, then we write~$N_G(v)$ and~$N_G[v]$ for the respective neighbourhoods. 
Slightly abusing notation, we may drop the subscript if the ambient graph is clear from context. 

For a positive integer~$d$, a set ${X \subseteq V(G)}$ is \emph{$d$-dominated} if there exists a set~${D \subseteq V(G)}$ with~${\abs{D} \leq d}$ such that~${X \subseteq N_G[D]}$. 

The \emph{average degree of~$G$} is the quantity~${2\abs{E(G)}/\abs{V(G)}}$. 
The \emph{independence number of~$G$}, denoted by~$\alpha(G)$, is the maximum size of an independent set in~$V(G)$. 

A graph~$H$ is a \emph{subdivision} of a graph~$G$ if~$H$ can be obtained from~$G$ by subdividing edges repeatedly, that is, removing an edge~$uv$ and adding new edges~$uw$ and $vw$ for some new vertex~$w$ not in~$G$. 

For a positive integer~$t$, let~$K_{t}$ denote the complete graph with~$t$ vertices. 
For positive integers~$s$ and~$t$, let~$K_{s,t}$ denote the complete bipartite graph with~$s$ vertices on one side and~$t$ vertices on the other side of the bipartition. 

\subsection{Tree-decompositions}

A \emph{tree-decomposition} of a graph~$G$ is a tuple~$(T,\beta)$ where~$T$ is a tree and~$\beta$ a map that maps every node~${t \in V(T)}$ to a subset of~$V(G)$ (called the \emph{bags} of the tree-decomposition) such that, for each edge~${uv \in E(G)}$, there is a node~${t \in V(T)}$ with~${u,v \in \beta(t)}$ and, for each~${v \in V(G)}$, the set of nodes~${t \in V(T)}$ for which~${v \in \beta(t)}$ induces a non-empty subtree of~$T$. 

Let~$G$ be a graph and let~$(T,\beta)$ be a tree-decomposition of~$G$. 
The \emph{width} of~$(T,\beta)$ is the maximum size of a bag minus one. 
The \emph{treewidth} of~$G$, denoted by~$\tw(G)$, is the minimum integer~$k$ such that~$G$ admits a tree-decomposition of width~$k$. 

For each edge~${st \in E(T)}$, the set~${\beta(s) \cap \beta(t)}$ is an \emph{adhesion set} of~$(T,\beta)$, and the \emph{adhesion} of~$(T,\beta)$ is the maximum size of an adhesion set. 

The \emph{independence number} of~$(T,\beta)$ is the maximum independence number of the graphs induced by the bags, that is, ${\max \{ \alpha(G[\beta(t)]) \colon t \in V(T) \}}$. 
The \emph{tree-independence number} of~$G$ is the minimum integer~$k$ such that~$G$ admits a tree-decomposition of independence number~$k$. 


\subsection{Clique minors}

We will need the following fundamental result of Duchet and Meyniel stating that every large graph either has a large independent set or a large clique minor. 

\begin{theorem}[Duchet, Meyniel \cite{DM82}]
    \label{thm:duchetmeyniel}
    For every graph~$G$, 
    \[
        (2 \cdot \alpha(G)-1) \cdot \had(G) \geq \lvert V(G) \rvert.
    \]
\end{theorem}

As a corollary, we immediately obtain \zcref{thm:atw}. 

\getkeytheorem{thmatw}

\begin{proof}
    Let~$(T,\beta)$ be a tree-decomposition of a graph~${G \in \mathcal{G}}$ of independence number at most~$c$. 
    Then, by \zcref{thm:duchetmeyniel}, ${\abs{\beta(t)} \leq (2c-1) \had(G[\beta(t)]) \leq (2c-1) \had(G)}$. 
    Hence, 
    \[
        \tw(G) \leq (2c-1) \had(G) - 1
    \]
\end{proof}

We will also need the following classical result of Kostochka~\cite{Kostochka:1984:HadwigerAverageDegree} and Thomason~\cite{Thomason:1984:ExtremalContractions,Thomason:2001:CompleteMinors} (see also \cite{AlonKS:2023:CompleteMinors} for a short proof).

\begin{theorem}\label{thm:kostochka-thomason}
    There exists an absolute constant~$a_{\ref{thm:kostochka-thomason}} > 0$ such that for every integer~${h \geq 2}$, 
    every graph~$G$ with~${\had(G) < h}$ has average degree 
    at most~${a_{\ref{thm:kostochka-thomason}} h \sqrt{\log h}}$. 
\end{theorem}

\subsection{Paths, linkages, and blocks}

Let~$G$ be a graph. 

A path~$P$ is an \emph{$(u,v)$-path} if both~$u$ and~$v$ are endvertices of~$P$. 
We call a set~$\mathcal{P}$ of pairwise internally vertex-disjoint $(u,v)$-paths an \emph{$(u,v)$-linkage}. 

Given sets~$A,B \subseteq V(G)$, we say that a path~$P$ is an \emph{$(A,B)$-path} if there are vertices~${a \in A}$ and~${b \in B}$ such that~$P$ is an~$(a,b)$-path whose internal vertices are disjoint from~${A \cup B}$.
An \emph{$(A,B)$-linkage} is a family of pairwise vertex-disjoint $(A,B)$-paths and a \emph{$(v,B)$-linkage} is a family of $(\{v\},B)$-paths which are pairwise vertex-disjoint except for~$v$. 
We call a set~${S \subseteq V(G)}$ an~\emph{$(A,B)$-separator} if~$S$ meets every~$(A,B)$-path. 

A set~$X$ is \emph{$k$-inseparable} if for every set~$S$ with~${\abs{S} < k}$, $X$ does not intersect two distinct components of~$G-S$. 
By Menger's Theorem, this means that between any two non-adjacent vertices~$u,v$ in~$X$ there is a $(u,v)$-linkage of size~$k$. 
A \emph{$k$-block} is a maximal $k$-inseparable set (of size at least~$k$). 
Slightly abusing notation, we say a \emph{$(k,s)$-block} is a $k$-inseparable set of size~$s$. 

From Menger's Theorem we can also directly observe the following. 
\begin{observation}
    \label{obs:blocks}
    If~$X$ is a $(k,2k)$-block in a graph~$G$ and~$(A,B)$ is a bipartition of~$X$ with~${\abs{A} = \abs{B} = k}$, then there is an~$(A,B)$-linkage of size~$k$ in~$G$. 
\end{observation}

\begin{proof}
    For a contradiction suppose there is no such linkage. 
    Then, by Menger's Theorem, there is a set~$S$ of size less than~$k$ such that no component of~${G - S}$ contains vertices from both~$A$ and~$B$. 
    But then there are~${u \in A}$ and~${v \in B}$ in different components of~${G-S}$, contradicting that~$X$ is $k$-inseparable. 
\end{proof}


A $(u,v)$-linkage between two nonadjacent vertices~$u$ and~$v$ is \emph{induced} if the paths in~$\mathcal{P}$ are induced and pairwise internally anti-complete. 


\begin{lemma}[{see also \cite[Corollary~11]{CampbellDDFGHHWWY:2024}}]
    \label{lem:inducedlinkage}
    Let~$h,k$ be positive integers, let~$G$ be a graph, let $u,v \in V(G)$ be nonadjacent, and let~$\mathcal{P}$ be a $(u,v)$-linkage in~$G$ of size at least~$2hk$.
    If~$\had(G) \leq h$, then~$G$ contains an induced $(u,v)$-linkage of size at least~$k$. 
\end{lemma}

\begin{proof}
    We define an auxiliary graph~$H$ with vertex set~$\mathcal{P}$ where~$P,Q$ are adjacent if and only if there exists an edge between an internal vertex of~$P$ and an internal vertex of~$Q$. 
    The auxiliary graph~$H$ is a minor of~$G$, so~${\had(H) \leq h}$. 
    Applying \zcref{thm:duchetmeyniel} to~$H$ yields a set~${\mathcal{Q} \subseteq \mathcal{P}}$ with~${\abs{\mathcal{Q}} \geq k}$ that is independent in~$H$. 
    By possibly shortening each of these paths to be induced, the resulting set~$\mathcal{R}$ is an induced $(u,v)$-linkage. 
    Note that since~$u$ and~$v$ are nonadjacent, each path in~$\mathcal{R}$ will contain an internal vertex from the original path in~$\mathcal{Q}$, so no two distinct paths in~$\mathcal{Q}$ yield the same path in~$\mathcal{R}$. 
\end{proof}

A \emph{strong $(k,s)$-block} is a $(k,s)$-block~$X$ such that $X$ is independent in~$G$ and between each pair~$u,v$ of vertices in~$X$ there is an induced $(u,v)$-linkage of size~$k$ consisting of induced paths. 

Chudnovsky et al.~\cite{ChudnovskyCFL:2025} showed that graphs with a large strong block contain either any planar graph or large complete bipartite graphs as induced minors. 
Their result is stated for hexagonal grids rather than arbitrary planar graphs. 
Every planar graph is an induced minor of a sufficiently large grid (see, e.g.,~\cite[Theorem~12]{CampbellDDFGHHWWY:2024}). 
Since every grid is an induced minor of a sufficiently large hexagonal grid, we may state their result as follows.

\begin{theorem}[{\cite[Theorem 5.1]{ChudnovskyCFL:2025}}]
    \label{thm:wallsKttslim}
    For every planar graph~$H$ and every positive integer~$t$, there are positive integers~$k$ and~$s$ such that the class $\Free{im}{H} \cap \Free{im}{K_{t,t}}$ contains no strong $(k,s)$-block. 
\end{theorem}

\smallskip

Given a positive integer~$k$, we say that a tree-decomposition~$(T,\beta)$ of a graph~$G$ is \emph{$k$-lean} if it has adhesion less than~$k$, and for any~${s, t \in V(T)}$ and~${A \subseteq \beta(s)}$, ${B \subseteq \beta(t)}$ with ${\ell \coloneqq \abs{A} = \abs{B} < k}$, there 
either is an $(A,B)$-linkage of size~${\ell}$ in~$G$, 
or there is an edge~${s't' \in E(T)}$ on the unique $(s,t)$-path in~$T$ such that the adhesion set~${\beta(s') \cap \beta(t')}$ has size less than~$\ell$. 
Since the adhesion of~$(T,\beta)$ is less than~$k$, each $k$-inseparable set of~$G$ is contained in some bag. 

In \cite{CarmesinDHH:2014:kblocks}, Carmesin et al.~noted that the proof given by Bellenbaum and Diestel~\cite{BellenbaumDiestel:2002} about the existence of lean tree-decompositions witnessing the treewidth of a finite graph can be adapted to prove the existence of a $k$-lean tree-decomposition of that graph.

\begin{theorem}[{\cite[Theorem 2.3]{CarmesinDHH:2014:kblocks}}]
    \label{thm:k-lean}
    Every graph admits a $k$-lean tree-decomposition for any positive integer~$k$. 
\end{theorem}

\subsection{Balanced separators and treewidth}

A \emph{weighting} of~$G$ is a map~${\mathsf{w} \colon V(G) \to \mathbb{R}_{\geq 0}}$ from~$V(G)$ to the non-negative reals. 
We call a weighting~${\mathsf{w} \colon V(G) \to \{ 0,1 \}}$ \emph{binary}. 
For a set~${S \subseteq V(G)}$, we write~$\mathsf{w}(S)$ for the sum~${\sum_{v \in S} \mathsf{w}(v)}$, and for a subgraph~${H \subseteq G}$, we write~$\mathsf{w}(H)$ for~$\mathsf{w}(V(H))$. 

For a weighting~$\mathsf{w}$ of~$G$ and~${c \in \halftoone}$, 
we call a subgraph~$H$ of~$G$ \emph{heavy} (with respect to~$\mathsf{w}$ and~$c$) if~${\mathsf{w}(V(H)) > c \cdot \mathsf{w}(V(G))}$. 
Clearly, for any set~${S \subseteq V(G)}$, the graph~${G - S}$ has at most one heavy component. 
We call a set~${S \subseteq V(G)}$ a \emph{$(\mathsf{w},c)$-balanced separator} if no component of~${G-S}$ is heavy. 

If~$S$ is a $d$-dominated $(\mathsf{w},c)$-balanced separator for a positive integer~$d$ and~${S \subseteq N_G[D]}$ for some set~${D \subseteq V(G)}$ with~${\abs{D} \leq d}$, then so is~${N_G[D]}$, so we often assume without loss of generality that our balanced separators have the form~$N_G[D]$. 

Graphs of bounded treewidth have such weighted balanced separators of bounded size. 
Indeed, given any tree-decomposition we can find a suitable bag which is a balanced separator.

\begin{lemma}[{\cite[Lemma 7.19]{CyganFKLMPPS:2015:ParameterizedAlgorithms}}]
    \label{lem:balancedbag}
    Every tree-decomposition~$(T,\beta)$ of a graph~$G$ with a weighting~$\mathsf{w}$ has a bag that is a $(\mathsf{w},\nicefrac{1}{2})$-balanced separator. 
\end{lemma}

On the other hand, the existence of such weighted balanced separators of bounded size provides a characterisation of graphs of bounded treewidth.
We use the following result of Harvey and Wood~\cite[Lemma~8]{HarveyWood:2017:TWtied}, formulated as below by Abrishami et~al.~\cite{AbrishamiCDHRSV:2024:ISTD2}. 

\begin{theorem}[{\cite[Lemma~8]{HarveyWood:2017:TWtied}; see also \cite[Lemma 1.13]{AbrishamiCDHRSV:2024:ISTD2}}]
    \label{thm:balancedseps}
    Let $G$ be a graph, let ${c \in \halftoone}$, and let $k$ be a positive integer. 
    If $G$ has a $(\mathsf{w}, c)$-balanced separator of size at most $k$ for every binary weighting~$\mathsf{w}$ on~$G$, then 
    \[
        \tw(G) \leq (1-c)^{-1}k\,.
    \]
\end{theorem}

\smallskip

Recall \zcref{conj:balancedseps}, stating that excluding any planar graph~$H$ as an induced minor yields a graph class with dominated binary-weighted balanced separators. 
This conjecture has been verified for several graphs~$H$, which we summarise as follows. 

\begin{theorem}
    \label{thm:dombalancedseps}
    For each of the following planar graphs~$H$, there is a positive integer~$d$ 
    such that the class~$\Free{im}{H}$ has $d$-dominated binary-weighted $\nicefrac{1}{2}$-balanced separators:
    \begin{itemize}
        \item $P_t$ for all positive integers~$t\geq 2$ with~${d = t-1}$ \textup{\cite{Gyarfas:1985:GyarfasPath}} \textup{(}see also~\textup{\cite{ChudnovskyPPT:2024:ApproxMWIS});} 
        \item $C_t$ for all positive integers~$t\geq3$ with~${d = t}$ \textup{\cite{ChudnovskyPPT:2024:ApproxMWIS} (}see also~\textup{\cite{ChudnovskyGKM:2025:WheelSeparators});}
        \item $F_t$ for all positive integers~$t\geq3$ with~${d = t}$ \textup{\cite{ChudnovskyGKM:2025:WheelSeparators};}
        \item $W_t$ for all positive integers~$t\geq3$ with~${d = (t-1)^2}$ \textup{\cite{ChudnovskyGKM:2025:WheelSeparators}.}
    \end{itemize}
\end{theorem}

Moreover, other natural classes of graphs have been shown to have dominated binary-weighted balanced separators, as we explain next.

A \emph{theta} is a subdivision of~$K_{2,3}$. 
A \emph{prism} is the line graph of a subdivision of~$K_{2,3}$. An \emph{even wheel} is a cycle of length at least~$4$ together with a vertex with an even number at least~$4$ of neighbours on the cycle. 
Note that while every $\ell$-wheel for an even integer~$\ell$ is an even wheel, not every even wheel is an $\ell$-wheel, as the cycle in an even wheel is allowed to contain more than~$\ell$ vertices. 
A graph is \emph{$C_4$-free odd-signable} if it excludes thetas, prisms, even wheels, and~$C_4$ as induced subgraphs. 

Recall that an \emph{even-hole-free graph} is a graph not containing any even cycle as an induced subgraph. 
It is easy to see that even-hole-free graphs do not contain any theta, prism, or even wheel as an induced subgraph; see e.g.~\cite{Aboulker:2013:Thesis}. 
Hence, every even-hole-free graph is $C_4$-free odd-signable.

\begin{theorem}[{\cite[Theorem 1.4]{ChudnovskyGHLS:2026:EvenHoleFree}}]
    \label{thm:dombalancedseps2}
    There is a positive integer~$d$ such that the class of $C_4$-free odd-signable graphs has $d$-dominated binary-weighted $\nicefrac{1}{2}$-balanced separators. In particular, the class of even-hole-free graphs has $d$-dominated binary-weighted $\nicefrac{1}{2}$-balanced separators. 
\end{theorem}

\section{The main method}
\label{sec:method}

The main method for proving our results follows the strategy we describe next. 
Suppose for now we have a method that, given a large $(k,s)$-block in our graph, allows us to construct a large clique minor (we will investigate this in the next section). 
Then, building upon a result of Chudnovsky, Gartland, Hajebi, Lokshtanov, and Spirkl~\cite{ChudnovskyGHLS:2026:EvenHoleFree}, we can find a block in the graph so that every balanced separator consisting of this block and the closed neighbourhood of a small set can be turned into a balanced separator consisting of the block and a set of small size. 
As we assume that we do not have large blocks, we can use \zcref{thm:balancedseps} to find our $(\tw,\had)$-bounding function. 

The result of Chudnovsky et al.~\cite[Theorem 9.2]{ChudnovskyGHLS:2026:EvenHoleFree} is stated in a different setting (with additional assumptions and a stronger conclusion), but their proof yields the following variant of their theorem. 
We include a proof for completeness. 

\begin{theorem}
    \label{thm:dom-to-small+block}
    Let~$k$ be a positive integer, let~$G$ be a graph, let~${c \in \halftoone}$, and let~$\mathsf{w}$ be a weighting of~$G$. 
    Then there exists a ${k}$-inseparable set~${K \subseteq V(G)}$ such that for every~${X \subseteq V(G) \setminus K}$ such that~$N[X] \cup K$ is a $(\mathsf{w},c)$-balanced separator there exists a set~$Y$ of size at most~$3 k \lvert X \rvert$ such that~$Y \cup K$ is a $(\mathsf{w},c)$-balanced separator. 
\end{theorem}

\begin{proof}
    The zero-weight case is immediate. 
    By \zcref{thm:k-lean}, $G$ admits a $(3k+1)$-lean tree-decomposition~${(T,\beta)}$. 
    By \zcref{lem:balancedbag}, this decomposition has a bag~$B$ which is a $(\mathsf w,\nicefrac{1}{2})$-balanced separator. 
    Let $K$ be the set of all vertices ${v \in V(G)}$ admitting a~$(v,B \setminus \{v\})$-linkage of size~$3k$. 
    
    For distinct ${u,v \in K}$, we can therefore find a set of~$3k$ $(u,v)$-walks by concatenating the respective linkages to~$B$ together with a linkage between their respective endpoints in~$B$ such that each vertex~${w \in V(G) \setminus \{u,v\}}$ is contained in at most~$3$ of the walks. 
    Now, any set ${S\subseteq V(G)\setminus\{u,v\}}$ with ${\abs{S} < k}$ meets at most~${3\abs{S} < 3k}$ of the walks, and hence, $K$ is $k$-inseparable. 

    For~${v \notin K}$, Menger's Theorem guarantees the existence of a set ${S_v \subseteq V(G)\setminus\{v\}}$ of size less than~$3k$ meeting every path from~$v$ to~${B\setminus\{v\}}$. 
    Now let~${X \subseteq V(G) \setminus K}$ be such that~${N[X] \cup K}$ is a $(\mathsf{w},c)$-balanced separator. 
    Let
    \[
        Y \coloneqq \bigcup_{v\in X}(S_v\cup\{v\}),
    \]
    and observe that~${\abs{Y} \leq 3k \abs{X}}$. 
    We claim that~${Y \cup K}$ is the desired $(\mathsf{w},c)$-balanced separator. 
    Suppose for a contradiction that a component~$C$ of~${G - (Y \cup K)}$ is heavy. 
    Then~${V(C) \cap B \neq \emptyset}$, since~${c \geq \nicefrac{1}{2}}$ and~$B$ is a $(\mathsf{w},\nicefrac{1}{2})$-balanced separator. 
    If some vertex~${v \in X}$ has a neighbour in~$C$, then there is a $(v,B\setminus\{v\})$-path avoiding~$S_v$, a contradiction. 
    Thus $C$ is contained in a component of~${G-(N[X] \cup K)}$, contradicting that~${N[X] \cup K}$ is a $(\mathsf{w},c)$-balanced separator. 
\end{proof}

With this, we can summarise our main method as follows. 

\begin{theorem}
    [store=mainmethod]
    \label{thm:mainmethod}
    Let~$\mathcal{C}$ be a hereditary class of graphs. 
    Suppose there is a~${c \in \halftoone}$ and there are functions~$d,b,s \colon \mathbb{N} \to \mathbb{N}$, with $d$ nondecreasing, such that every graph~$G$ in~$\mathcal{C}$
    \begin{enumerate}
        \item has a $d(\had(G))$-dominated $(\mathsf{w},c)$-balanced separator~$S$ for every binary weighting~$\mathsf{w}$ of~$G$ and
        \item contains no $(b(\had(G)),s(\had(G)))$-block. 
    \end{enumerate}
    Then the function~$f \colon \mathbb{N} \to \mathbb{R}_{\geq 0}$ defined by 
    \[
        f(h) \coloneqq (1-c)^{-1} ( 3 \cdot b(h) \cdot d(h) + s(h) )
    \]
    is a $(\tw,\had)$-bounding function. 
\end{theorem}

\begin{proof}
    Let~${G \in \mathcal{C}}$, $h \coloneqq \had(G)$, ${k \coloneqq b(h)}$, ${s \coloneqq s(h)}$, ${d \coloneqq d(h)}$, and let~$\mathsf{w}$ be a binary weighting of~$G$.
    
    We apply \zcref{thm:dom-to-small+block} to get a $k$-inseparable set~$K$ with $\abs{K} < s$. 

    Let~$S$ be a $(\mathsf{w},c)$-balanced separator of~${G-K}$ dominated by a set~$X$ of size at most~$d$.
    Such a separator exists since~$\mathcal{C}$ is hereditary and the function~$d$ is monotone, using the restriction of $\mathsf{w}$ to $G-K$. 
    Thus $N_G[X] \cup K$ is a $(\mathsf{w},c)$-balanced separator of~$G$. 
    Then by \zcref{thm:dom-to-small+block} there exists a $(\mathsf{w},c)$-balanced separator of size less than~${3 k d + s}$. 
    Since this holds for every binary weighting, by \zcref{thm:balancedseps}, 
    \[
        \tw(G) \leq (1-c)^{-1}(3kd+s) = f(\had(G)),
    \]
    as desired. 
\end{proof}

\section{From blocks to clique minors}
\label{sec:blocks2cliqueminors}

To be able to use our main method, we want to argue that if we have a $(k,s)$-block for sufficiently large~$k$ and~$s$, then we can find a clique minor in graphs from the classes we are interested in. 
In this section, we will do exactly that for several classes of graphs. 

Let us start with graphs that do not contain induced subdivisions of~$K_{2,t}$. 

\begin{lemma}
    \label{lem:K2t-sub-free}
    Let~$t, h$ be positive integers. 
    Then every graph~$G$ that contains a $(2th,h)$-block either contains
    \begin{enumerate}
        \item a subdivision of~$K_{2,t}$ as an induced subgraph, or
        \item a $K_h$-minor. 
    \end{enumerate}
\end{lemma}

\begin{proof}
    Let~$X$ be a $(2th,h)$-block. 
    Assume that $\had(G) < h$. 
    Then~$X$ is not a clique, so we may assume there are two non-adjacent vertices~${x,y \in X}$ and a $(x,y)$-linkage~$\mathcal{P}$ of size $2th$, by Menger's Theorem. 
    By \zcref{lem:inducedlinkage} there is an induced~$(x,y)$-linkage~$\mathcal{Q}$ of size~$t$. 
    But then, $\bigcup \mathcal{Q}$ is an induced subgraph of~$G$ which is a subdivision of~$K_{2,t}$. 
\end{proof}

The next lemma provides a tool for finding, in the absence of a large clique minor, a strong block within a block. 

\begin{lemma}
    \label{lem:strongblock}
    Let~$h,\ell,k,s$ be positive integers. 
    Then every graph~$G$ that contains a $(2kh^{\ell},2sh^\ell)$-block either contains
    \begin{enumerate}
        \item a strong $(kh^{\ell-1},sh^{\ell-1})$-block, or
        \item a $K_h$-minor. 
    \end{enumerate}
\end{lemma}

\begin{proof}
    Let~$X$ be the $(2kh^{\ell},2sh^{\ell})$-block and suppose ${\had(G) < h}$. 
    As before, for each pair~${x,y \in X}$ of non-adjacent vertices there is a $(x,y)$-linkage of size~${2kh^\ell}$ and hence an induced $(x,y)$-linkage of size~${kh^{\ell-1}}$ by \zcref{lem:inducedlinkage}. 
    By again applying \zcref{thm:duchetmeyniel} to~$G[X]$, we get that
    \[
        \alpha(G[X]) \geq \frac{\abs{X}}{2h} \geq sh^{\ell-1}. 
    \]
    Now, any independent set~${I \subseteq X}$ with~${\abs{I} = sh^{\ell-1}}$ is a strong $(kh^{\ell-1},sh^{\ell-1})$-block. 
\end{proof}

Using \zcref{thm:wallsKttslim}, this now allows us to find, in a graph with a large block, either a planar graph as an induced minor, or a large biclique as an induced minor, or a large clique minor. 

\begin{lemma}
    \label{lem:H-Ktt-im-free}
    Let~$t$ be a positive integer and let~$H$ be a planar graph. 
    Then there are positive integers~$k, s$ such that, for every positive integer~$h$, every graph~$G$ that contains a $(kh,sh)$-block either contains
    \begin{enumerate}
        \item $H$ as an induced minor, 
        \item $K_{t,t}$ as an induced minor, or
        \item $K_h$ as a minor. 
    \end{enumerate}
\end{lemma}

\begin{proof}
    Let~${\mathcal{C} \coloneqq \Free{im}{H} \cap \Free{im}{K_{t,t}}}$. 
    By \zcref{thm:wallsKttslim},
    there exist positive integers~$p,q$ such that no graph in~$\mathcal{C}$ contains any strong~$(p,q)$-blocks. 
    Let~${k \coloneqq 2p}$ and~${s \coloneqq 2q}$. 
    So assume~${G \in \mathcal{C}}$ contains a $(kh,sh)$-block. 
    Then, by \zcref{lem:strongblock}, $G$ contains a $K_h$-minor. 
\end{proof}

For a positive integer~$d$ and a~${c \in \halftoone}$, let us call a set~${X \subseteq V(G)}$ \emph{$(c,d)$-separable} if for every induced subgraph~${H \subseteq G}$ and binary weighting~$\mathsf{w}$ of~$H$ whose support is a subset of~$X$ there is a $(\mathsf{w},c)$-balanced separator in~$H$ dominated by at most~$d$ vertices. 
Note that if~$G$ is in a hereditary class with $d$-dominated binary-weighted $c$-balanced separators, then every~${X \subseteq V(G)}$ is $(c,d)$-separable. 

\begin{theorem}
    [store=balancedblock]
    \label{thm:balancedblock}
    For every positive integer~$d$ and every~${c \in \halftoone}$ there exists a positive integer~$a_{\ref{thm:balancedblock}}(c,d)$ such that for every integer~${h \geq 2}$ the following statement holds. 
    Let
    \[
        k\coloneqq a_{\ref{thm:balancedblock}}(c,d)\, h\left\lceil\log h\right\rceil.
    \] 
    Then every graph~$G$ that contains a $(k,2k)$-block that is $(c,d)$-separable contains a $K_h$-minor. 
\end{theorem}

We postpone the proof of this theorem to \zcref{sec:balancedblock}. 

\section{Combining the results}
\label{sec:applications}

In this section, we prove our main results by combining the results of \zcref{sec:blocks2cliqueminors} with our main method, \zcref{thm:mainmethod}. 

\getkeytheorem{thmquasi}

\begin{proof}
    Let $a \coloneqq a_{\ref{thm:balancedblock}}(c,d)$ be the constant in \zcref{thm:balancedblock}, and set
    \[
        b(h) \coloneqq a(h+1)\left\lceil\log(h+1)\right\rceil,
        \qquad 
        s(h) \coloneqq 2b(h).
    \]
    Every subset of a graph in~$\mathcal{G}$ is $(c,d)$-separable since~$\mathcal{G}$ is hereditary and has $d$-dominated binary-weighted $c$-balanced separators.
    Hence a graph ${G \in \mathcal{G}}$ with ${\had(G) = h}$ has no $(b(h),s(h))$-block, since such a block would force a $K_{h+1}$-minor by \zcref{thm:balancedblock}. 
    Now \zcref{thm:mainmethod} gives
    \[
        \tw(G)
        \leq
        \mbox{$\frac{3d+2}{1-c}$}\,b(h)
        \in O(h \log h),
    \]
    proving the asserted asymptotic bound.
\end{proof}

\zcref{cor:main+conj} now is an immediate consequence of \zcref{thm:main,thm:dombalancedseps}. 

Let us now verify the linear bounds we claimed in the introduction. 

\getkeytheorem{thmevenholefree}

\begin{proof}
    By \zcref{thm:dombalancedseps2}, the class~$\mathcal{C}$ of $C_4$-free odd-signable graphs has $d$-dominated binary-weighted $\nicefrac{1}{2}$-balanced separators for some positive integer~$d$. 
    By definition of~$\mathcal{C}$, no graph in~$\mathcal{C}$ contains any subdivision of~$K_{2,3}$ as an induced subgraph. 
    Set ${b(h) \coloneqq 6(h+1)}$ and~${s(h) \coloneqq h+1}$. 
    Now, no~${G \in \mathcal{C}}$ with~$\had(G) = h$ has a $(b(h),s(h))$-block by \zcref{lem:K2t-sub-free}. 
    \zcref{thm:mainmethod} yields
    \[
        \tw(G) \leq (36d+2)\,(h+1).
    \]
    And as we remarked in the introduction, it is easy to see that every even-hole-free graph is $C_4$-free odd-signable; see e.g.~\cite{Aboulker:2013:Thesis}. 
\end{proof}

\getkeytheorem{thmHandKtt}

\begin{proof}
    Let $k_0,s_0$ be the constants of \zcref{lem:H-Ktt-im-free} and take $b(h) \coloneqq k_0(h+1)$ and $s(h) \coloneqq s_0(h+1)$.
    \zcref{thm:mainmethod} yields
    \[
        \tw(G) \leq \mbox{$\frac{3dk_0 + s_0}{1-c}$}\,(h+1) \qedhere
    \]
\end{proof}

\section{\texorpdfstring{Proof of \zcref{thm:balancedblock}}{Proof of Theorem 4.4}}
\label{sec:balancedblock}

For the proof of \zcref{thm:main} to be complete, we have to provide a proof of \zcref{thm:balancedblock}. 
Towards this goal, let us first extract some structure from a graph that contains a large~$(c,d)$-separable block. 

\begin{lemma}
    \label{lem:separator-cut-family}
    Let~${c \in \halftoone}$, let~${d \geq 1}$ and~${k \geq 16d/(1-c)}$ be integers, 
    and let~$X$ be a $(c,d)$-separable $(k,2k)$-block in~$G$. 
    Let~${t \coloneqq \left\lfloor \frac{(1-c) k}{8d} \right\rfloor}$. 
    Then there are 
    \begin{itemize}
        \item $t$ pairwise disjoint nonempty sets~${D_1,\dots,D_t}$, each of size at most~$d$, 
    \end{itemize}
    and with~${Z \coloneqq \bigcup_{i\in[t]} D_i}$ and~${q \coloneqq \abs{Z}}$, 
    there is
    \begin{itemize}
        \item a set~${M \subseteq X \setminus Z}$ of size~${m \coloneqq 2k-2q}$ that is an ${(m/2,m)}$-block in~${G-Z}$,
    \end{itemize}
    such that for every~${i \in [t]}$, 
    \begin{enumerate}
        \item either~$D_i$ has at least~${(1-c) k/4}$ neighbours in~$M$, or
        \item there is a bipartition~${\{A_i,B_i\}}$ of~$M$ with
            \[
                \min\{\abs{A_i},\abs{B_i}\} 
                \geq (1-c) k/2,
            \]
            such that~${N_G(D_i) \setminus Z}$ is an $(A_i,B_i)$-separator in~${G-Z}$. 
    \end{enumerate}

    Moreover, we have
    \begin{itemize}
        \item $\frac{(1-c) k}{16d} \leq t \leq q \leq \frac{(1-c) k}{8} \leq \frac{k}{16}$. 
    \end{itemize}
\end{lemma}

\begin{proof}    
    Let us first construct the family~$D_1, \dots, D_t$ recursively. 
    Let~${i \in [t]}$ and suppose~$D_j$ has already been chosen for each~${j < i}$. 
    Let~${Z_i \coloneqq \bigcup_{j<i} D_j}$. 
    We apply the $(c,d)$-separability of~$X$ in~${G-Z_i}$ for the binary weighting~$\mathsf{w}_i$ with support~${X \setminus Z_i}$ to obtain a set~${D_i \subseteq V(G) \setminus Z_i}$ with~${\abs{D_i} \leq d}$, such that
    \[
        S_i \coloneqq N_{G-Z_i}[D_i]
    \]
    is a $(\mathsf{w}_i,c)$-balanced separator of~${G - Z_i}$. 
    Since~${\abs{Z_i} \leq (i-1)d < td \leq \frac{(1-c)k}{8} < k}$, all vertices of~${X \setminus Z_i}$ lie in the same component of~${G-Z_i}$.
    Thus~$D_i$ cannot be empty. 
    
    We easily observe that~${{(1-c) k}/{(16d)} \leq t \leq {(1-c) k}/{(8d)}}$. 
    Now the sets~${D_1,\dots,D_t}$ are nonempty, pairwise disjoint, of size at most~$d$, and we have~${t \leq q \leq td \leq (1-c) k/8}$.

    Since~${\abs{X \setminus Z} \geq 2k-q}$, let~$M$ be any subset of~${X \setminus Z}$ of size~${m = 2k-2q}$. 
    If a set~${S \subseteq V(G-Z)}$ has size less than~${m/2 = k-q}$, then~${\abs{S \cup Z} < k}$.
    And since~$X$ is $k$-inseparable in~$G$, it follows that~$M$ is a $(m/2,m)$-block in~${G-Z}$. 

    Now fix an index~${i \in [t]}$. 
    Since~$S_i$ is a $(\mathsf{w}_i,c)$-balanced separator of~$G-Z_i$, by the definition of~$\mathsf{w}_i$, each component of~${(G-Z_i)-S_i}$ contains at most~${c\abs{X \setminus Z_i} \leq c2k}$ vertices of~$M$. 
    If~${\abs{S_i \cap M} \geq (1-c) k/4}$, the first conclusion holds since~${D_i \subseteq Z}$ is disjoint from~$M$. 
    Otherwise, we select components of~${(G-Z_i)-S_i}$ in decreasing order of their number of vertices in~$M$, until their union first contains at least~${(1-c) k/2}$ such vertices in~$M$. 
    This is possible, since~$2q \leq (1-c)k/4$ and hence
    \[
        \abs{M \setminus S_i} 
        > (2k - 2q) - (1-c) k/4
        \geq 2k - (1-c) k/2 
        > (1-c) k/2.
    \]
    Let~$A_i$ be the intersection of this union with~$M$.
    If the first component suffices, then~${\abs{A_i} \leq c2k}$. 
    Otherwise, every selected component has fewer than~${(1-c) k/2}$ vertices of~$M$, so~${\abs{A_i} < (1-c) k \leq c2k}$.
    Hence, 
    \[
        (1-c) k/2 
        \leq \abs{A_i} 
        \leq c2k.
    \]
    Then~${B_i \coloneqq M \setminus A_i}$ satisfies, again using $2q \leq (1-c)k/4$, 
    \begin{align*}
        \abs{B_i}
        &= \abs{M} - \abs{A_i}
        \geq (2k - 2q) - c2k
        \geq 2k ( 1 - (1-c)/8 - c ) \\
        &= 2k ( (1 - c) - (1-c)/8 )
        \geq 7 (1-c) k / 4
        \geq (1-c) k / 2, 
    \end{align*}
    as required. 
    By construction, any $(A_i,B_i)$-path in~${G-Z}$ meets~$S_i$. 
    As the path avoids~${D_i \subseteq Z}$, it has a vertex adjacent to~$D_i$. 
\end{proof}

We use the following form of the Arora--Rao--Vazirani geometric separation theorem. 

\begin{theorem}[{Arora, Rao, Vazirani~\cite[Theorem~1]{AroraRV:2009:ExpanderFlows}}]
    \label{thm:ARV-separation}
    For every~${\beta \in (0,\nicefrac{1}{2})}$ there are constants~${\delta,\eta > 0}$ with the following property. 
    Let~$n \geq 2$ be an integer and let~${v_1,\ldots,v_n}$ be unit vectors in $n$-dimensional Euclidean space such that for each~${i,j,k \in [n]}$, 
    \[
        \norm{v_i - v_j}^2 + \norm{v_j - v_k}^2 \geq \norm{v_i - v_k}^2
    \]
    and
    \[
        \sum_{i<j} \norm{v_i-v_j}^2
        \geq 4 \beta (1-\beta) n^2\,.
    \]
    Then there are disjoint index sets~${U,V \subseteq[n]}$, each of size at least~${\delta n}$, such that
    \[
        \norm{v_u - v_v}^2
        \geq \frac{\eta}{\sqrt{\log n}}
        \qquad\text{for every }u\in U,\ v\in V.
    \]
\end{theorem}

This theorem, together with our structural observations from \zcref{lem:separator-cut-family}, will allow us to find a minor of large average degree. 

\begin{lemma}
    \label{lem:ARV-dense-minor}
    For every~${c \in \halftoone}$ there is a constant~$a_{\ref{lem:ARV-dense-minor}}(c)>0$ satisfying the following property. 
    Let~${d \geq 1}$ and~${k \geq 16 d / (1-c)}$ be integers and let~$X$ be a $(c,d)$-separable $(k,2k)$-block in~$G$. 
    Then~$G$ has a minor with at most~$2k$ vertices and average degree at least~${{a_{\ref{lem:ARV-dense-minor}}(c)k}/{({d\sqrt{\log k}})}}$. 
\end{lemma}

\begin{proof}
    Let~${\beta \coloneqq (1-c)/4}$ and let~$\delta, \eta > 0$ be the constants given by \zcref{thm:ARV-separation} for this choice of~$\beta$. 
    Note that~$\delta,\eta$ depend only on~$c$. 
    We define
    \[
        a_{\ref{lem:ARV-dense-minor}}(c) 
        \coloneqq 
        \min\left\{
        \frac{(1-c)^2}{64},
        \frac{(1-c)\delta\eta}{64} 
        \right\} \, .
    \]

    We apply \zcref{lem:separator-cut-family} to find for~${t \coloneqq \left\lfloor \frac{(1-c) k}{8d} \right\rfloor}$, 
    \begin{itemize}
        \item $t$ pairwise disjoint nonempty sets~${D_1,\dots,D_t}$, each of size at most~$d$, 
    \end{itemize}
    and with~${Z \coloneqq \bigcup_{i\in[t]} D_i}$ and~${q \coloneqq \abs{Z}}$, 
    \begin{itemize}
        \item a set~${M \subseteq X \setminus Z}$ of size~${m \coloneqq 2k-2q}$ that is an ${(m/2,m)}$-block in~${G-Z}$, 
    \end{itemize}
    satisfying the properties stated in the lemma. 
    
    Suppose first that ${\abs{N_G(D_i) \cap M} \geq (1-c)k/4}$ for at least~${t/2}$ indices~${i \in [t]}$. 
    Then, since the sets~$D_1, \dots, D_t$ are pairwise disjoint, there are at least~${t (1-c) k / 8}$ distinct edges between~$Z$ and~$M$. 
    For the induced subgraph~$J \coloneqq G[Z \cup M]$, we have
    \[
        \abs{V(J)} = q + m = q + 2k - 2q 
        = 2k - q \leq 2k
        \,,
    \]
    and since~${\sqrt{\log k} \geq \sqrt{\log 32} \geq 2}$, the average degree of~$J$ is 
    \[
        \frac{2 \abs{E(J)}}{\abs{V(J)}}
        \geq \frac{2 \cdot (1-c)tk}{8 \cdot 2 k} 
        = \frac{(1-c) t}{8}
        \geq \frac{(1-c)^2 k}{8 \cdot 16 d}
        \geq \frac{a_{\ref{lem:ARV-dense-minor}}(c)k}{d\sqrt{\log{k}}}
        \,,
    \]
    so~$J$ is the desired minor. 

    Otherwise, let~${I \subseteq [t]}$ be the set of indices for which the second conclusion holds, that is, for each~${i \in I}$ there is a bipartition~${\{A_i,B_i\}}$ of~$M$ with~${\min\{ \abs{A_i}, \abs{B_i}\} \geq (1-c)k/2 \geq 2\beta k}$ such that~${S_i \coloneqq N_G(D_i) \setminus Z}$ is an~$(A_i,B_i)$-separator in~${G-Z}$. 
    Let~${s \coloneqq \abs{I} \geq t/2}$.
    
    For each~${x \in M}$, define the unit vector
    \[
        v_x \coloneqq \frac{1}{\sqrt{s}} \bigl( \sigma_i(x) \bigr)_{i \in I},
        \qquad
        \sigma_i(x) \coloneqq
        \begin{cases}
             1,  & x \in A_i,\\
            -1,  & x \notin A_i.
        \end{cases}
    \]
    For each~${\{x,y\} \subseteq M}$, let~$I_{\{x,y\}} \subseteq I$ denote the set of indices for which~$x$ and~$y$ are on different sides of the bipartition~$\{A_i,B_i\}$, that is, ${I_{\{x,y\}} \coloneqq \{ i \in I \colon \abs{A_i \cap \{x,y\}} = 1 \}}$. 
    Then, 
    \begin{equation}\label{eq:separator-cut-distance}
        \norm{v_x-v_y}^2
        = \mbox{$\frac{4}{s}$} \abs{ I_{\{x,y\}} }.
    \end{equation}
    Note that for distinct~${x,y,z \in M}$, every index~${i \in I_{\{x,z\}}}$ belongs to exactly one of~$I_{\{x,y\}}$ and~$I_{\{y,z\}}$. 
    Hence, $I_{\{x,z\}} \subseteq I_{\{x,y\}} \cup I_{\{y,z\}}$, 
    and so $\abs{I_{\{x,z\}}} \leq \abs{I_{\{x,y\}}} + \abs{I_{\{y,z\}}}$. 
    Together with \eqref{eq:separator-cut-distance}, this gives
    \[
        \norm{v_x-v_y}^2+\norm{v_y-v_z}^2
        \geq
        \norm{v_x-v_z}^2.
    \]
    For every~${i \in I}$, since ${\min\{ \abs{A_i}, \abs{B_i}\} \geq 2\beta k \geq \beta m}$ and ${\abs{A_i} = m - \abs{B_i}}$, we obtain that ${\beta m \leq \abs{A_i} \leq (1-\beta)m}$. 
    Moreover, since the function ${x \mapsto x(m-x)}$ is concave, its minimum on ${[\beta m, (1-\beta)m]}$ occurs at an endpoint. 
    Hence, \[{\abs{A_i}\abs{B_i} = \abs{A_i}(m-\abs{A_i}) \geq \beta(1-\beta)m^2}.\] 
    Consequently,  
    \begin{align*}
        \sum_{\{x,y\} \subseteq M} \norm{v_x-v_y}^2
        &= \mbox{$\frac{4}{s}$} \sum_{i \in I} (\abs{A_i} \cdot \abs{B_i})  
        \geq 4\beta(1-\beta)m^2. 
    \end{align*}
    Thus all hypotheses in the premise of \zcref{thm:ARV-separation} hold, so there exist disjoint sets~${U,V \subseteq M}$, each of size at least~$\delta m$, such that for every pair~${(x,y) \in U \times V}$, we have~$\norm{v_x - v_y}^2 \geq \eta/\sqrt{\log m}$. 
    Using \eqref{eq:separator-cut-distance}, we conclude that for every such pair, 
    \[
        \abs{I_{\{x,y\}}} \geq \frac{\eta s}{4\sqrt{\log m}}.
    \]

    Let~${\ell \coloneqq \min\{ \abs{U}, \abs{V} \}}$.
    Since~$M$ is $(m/2)$-inseparable in~$G-Z$ and~$\ell \leq m/2$, Menger's Theorem gives an $(U,V)$-linkage~$\mathcal{P}$ of size~${\ell \geq \delta m}$ in~${G-Z}$.
    Each path in~$\mathcal{P}$ with endpoints~${x \in U}$ and~${y \in V}$ crosses the bipartition~${\{A_i,B_i\}}$ for all~${i \in I_{\{x,y\}}}$ and consequently contains a vertex in~${S_i = N_G(D_i) \setminus Z}$. 

    Now we obtain a model of a minor~$J$ with vertex set~${\mathcal{P} \cup Z}$ using the branch sets ${\{ V(P) \colon P \in \mathcal{P}\} \cup \{ \{ z \} \colon z \in Z \}}$ where two vertices of~$J$ are adjacent if there is an edge between the corresponding branch sets. 
    Note that since each path in~$\mathcal{P}$ is a path in~${G-Z}$, the branch sets are indeed disjoint.  
    As discussed above, ${P \in \mathcal{P}}$ with endpoints~${x \in U}$ and~${y \in V}$ is adjacent in~$J$ to a vertex in~${D_i \subseteq Z}$ for all~${i \in I_{\{x,y\}}}$. 
    Since for these pairs~$\{x,y\}$ we have~$\abs{I_{\{x,y\}}} \geq \eta s / (4 \sqrt{ \log m})$,
    we obtain 
    \[
        \abs{V(J)} = q+\ell \leq k,
        \qquad
        \abs{E(J)} 
        \geq \frac{\ell\eta s}{4\sqrt{\log m}}. 
    \]
    Using~$s\geq t/2\geq (1-c) k/(32d)$, and~${\ell\geq\delta m}$, we get
    \[
        \frac{2|E(J)|}{|V(J)|}
        \geq\frac{\ell\eta s}{2k\sqrt{\log m}}
        \geq\frac{(1-c) \delta\eta m}{64d\sqrt{\log m}} 
        \geq \frac{a_{\ref{lem:ARV-dense-minor}}(c)k}{d\sqrt{\log{k}}}
        \,,
    \]
    where the last inequality holds since for~${x \geq 2}$, the function~${x \mapsto x/\sqrt{\log x}}$ is increasing and~${m \geq k}$. 
    
    Hence, in both cases, we have the desired minor~$J$. 
\end{proof}

We are finally ready to prove \zcref{thm:balancedblock}, which we restate for convenience. 

\getkeytheorem{balancedblock}

\begin{proof}
    Let~$a_{\ref{thm:kostochka-thomason}}$ be the constant from \zcref{thm:kostochka-thomason} and let~$a_{\ref{lem:ARV-dense-minor}}(c)$ be the constant from \zcref{lem:ARV-dense-minor}. 
    Choose an integer~$a\geq8d/(1-c)$ such that
    \[
        a>\frac{d a_{\ref{thm:kostochka-thomason}} }{a_{\ref{lem:ARV-dense-minor}}(c)}\sqrt{3+\log a}.
    \]
    Note that such an integer exists since
    ${x/\sqrt{3+\log x}\to\infty}$ as~${x \to \infty}$. 
    Fix~${h \geq 2}$ and set~${k \coloneqq ah \lceil\log h\rceil}$. 
    Then ${2k \geq 32d/(1-c)}$, so \zcref{lem:ARV-dense-minor} gives a minor~$J$ of average degree at least~$\frac{a_{\ref{lem:ARV-dense-minor}}(c)k}{d\sqrt{\log k}}$.
    Suppose for a contradiction that~$\had(G)<h$.
    Since~$\had(J)<h$, \zcref{thm:kostochka-thomason} gives
    \[
        k\leq\frac{d a_{\ref{thm:kostochka-thomason}}}{a_{\ref{lem:ARV-dense-minor}}(c)}
        h\sqrt{\log h}\sqrt{\log k}.
    \]
    As~$\lceil\log h\rceil\leq2h$ and~$\log h\geq1$, we have
    \begin{align*}
        \log k 
        &= \log a + \log h + \log \lceil\log h\rceil  
        \leq \log a + \log 2 + 2 \log h\\ 
        &= \log a+1+2\log h 
        \leq (3+\log a)\log h.
    \end{align*}
    Consequently,
    \[
        k\leq\frac{d a_{\ref{thm:kostochka-thomason}}}{a_{\ref{lem:ARV-dense-minor}}(c)}\sqrt{3+\log a}\,h\log h
        <ah\log h\leq k,
    \]
    a contradiction.
\end{proof}

\section*{Acknowledgements}

This work was initiated at the workshop \emph{Structural Graph Theory Workshop
on the Adriatic Coast}, SGTAC Workshop 2026, hosted by the University of Primorska. 
The authors would like to thank the organisers of the workshop as well as all participants for providing a productive research atmosphere. 
In particular, the authors thank Maria Chudnovsky for pointing out some results that turned out to be helpful for their approach. 

J.~Pascal Gollin was supported in part by the Slovenian Research and Innovation Agency (research project N1-0370).

\section*{AI declaration}

This work was assisted by ChatGPT-5.6 Sol and ChatGPT-6 Astra. 
Specifically, the models pointed out the connection to the geometric separation theorem of Arora, Rao, and Vazirani (\zcref{thm:ARV-separation}), and provided an initial proof of~\zcref{thm:balancedblock}, 
which was subsequently edited and improved by the authors. The models were also used for proofreading. 
The authors take full responsibility for the content of the paper.

\printbibliography

@unpublished{CampbellDDFGHHWWY:2024,
    AUTHOR = {Campbell, Rutger and Davies, James and Distel, Marc and Frederickson, Bryce and Gollin, J. Pascal and Hendrey, Kevin and Hickingbotham, Robert and Wiederrecht, Sebastian and Wood, David and Yepremyan, Liana},
    TITLE = {Treewidth, {H}adwiger number, and induced minors},
    YEAR = {2024},
    NOTE = {Preprint},
    EPRINT = {2410.19295},
    EPRINTTYPE = {arxiv},
}

@unpublished{ChudnovskyCFL:2025,
    AUTHOR = {Chudnovsky, Maria and Codsi, Julien and Fischer, David and Lokshtanov, Daniel},
    TITLE = {Induced minors and subpolynomial treewidth},
    YEAR = {2025},
    NOTE = {Preprint},
    EPRINT = {2512.18835},
    EPRINTTYPE = {arxiv},
}

@incollection {DM82,
    AUTHOR = {Duchet, Pierre and Meyniel, Henri},
     TITLE = {On {H}adwiger's number and the stability number},
 BOOKTITLE = {Graph theory ({C}ambridge, 1981)},
    SERIES = {North-Holland Math. Stud.},
    VOLUME = {62},
     PAGES = {71--73},
 PUBLISHER = {North-Holland, Amsterdam-New York},
      YEAR = {1982},
      ISBN = {0-444-86449-0},
   MRCLASS = {05C35 (05C15 05C75)},
  MRNUMBER = {671905},
       DOI = {10.1016/S0304-0208(08)73549-7},
}

@article {DvorakNorin:2019:BalancedSeps,
    AUTHOR = {Dvo\v r\'ak, Zden\v ek and Norin, Sergey},
     TITLE = {Treewidth of graphs with balanced separations},
   JOURNAL = {J. Combin. Theory Ser. B},
  FJOURNAL = {Journal of Combinatorial Theory. Series B},
    VOLUME = {137},
      YEAR = {2019},
     PAGES = {137--144},
      ISSN = {0095-8956,1096-0902},
   MRCLASS = {05C75 (05C05 05C70)},
  MRNUMBER = {3980088},
MRREVIEWER = {Yan\ Liu},
       DOI = {10.1016/j.jctb.2018.12.007},
       URL = {https://doi.org/10.1016/j.jctb.2018.12.007},
    EPRINT = {1408.3869},
EPRINTTYPE = {arxiv},
}

@article {HarveyWood:2017:TWtied,
    AUTHOR = {Harvey, Daniel J. and Wood, David R.},
     TITLE = {Parameters tied to treewidth},
   JOURNAL = {J. Graph Theory},
  FJOURNAL = {Journal of Graph Theory},
    VOLUME = {84},
      YEAR = {2017},
    NUMBER = {4},
     PAGES = {364--385},
      ISSN = {0364-9024,1097-0118},
   MRCLASS = {05C85 (05C75)},
  MRNUMBER = {3623383},
MRREVIEWER = {Jisu\ Jeong},
       DOI = {10.1002/jgt.22030},
       URL = {https://doi.org/10.1002/jgt.22030},
    EPRINT = {1312.3401},
EPRINTTYPE = {arxiv},
}

@article {AbrishamiCDHRSV:2024:ISTD2,
    AUTHOR = {Abrishami, Tara and Chudnovsky, Maria and Dibek, Cemil and
              Hajebi, Sepehr and Rzążewski, Pawe\l{} and Spirkl,
              Sophie and Vu\v skovi\'c, Kristina},
     TITLE = {Induced subgraphs and tree decompositions {II}. {T}oward walls
              and their line graphs in graphs of bounded degree},
   JOURNAL = {J. Combin. Theory Ser. B},
  FJOURNAL = {Journal of Combinatorial Theory. Series B},
    VOLUME = {164},
      YEAR = {2024},
     PAGES = {371--403},
      ISSN = {0095-8956,1096-0902},
   MRCLASS = {05C75 (05C76)},
  MRNUMBER = {4660624},
MRREVIEWER = {Andr\'e\ E.\ K\'ezdy},
       DOI = {10.1016/j.jctb.2023.10.005},
       URL = {https://doi.org/10.1016/j.jctb.2023.10.005},
    EPRINT = {2108.01162},
EPRINTTYPE = {arxiv},
}

@article {AroraRV:2009:ExpanderFlows,
    AUTHOR = {Arora, Sanjeev and Rao, Satish and Vazirani, Umesh},
     TITLE = {Expander flows, geometric embeddings and graph partitioning},
   JOURNAL = {J. ACM},
  FJOURNAL = {Journal of the ACM},
    VOLUME = {56},
      YEAR = {2009},
    NUMBER = {2},
     PAGES = {Art. 5, 37},
      ISSN = {0004-5411,1557-735X},
   MRCLASS = {68R10 (05C85 60D05 68W25)},
  MRNUMBER = {2535878},
       DOI = {10.1145/1502793.1502794},
       URL = {https://people.eecs.berkeley.edu/~vazirani/pubs/arvjacm.pdf},
    EPRINT = {},
EPRINTTYPE = {},
      NOTE = {},
}

@article {Kostochka:1984:HadwigerAverageDegree,
    AUTHOR = {Kostochka, A. V.},
     TITLE = {Lower bound of the {H}adwiger number of graphs by their
              average degree},
   JOURNAL = {Combinatorica},
  FJOURNAL = {Combinatorica. An International Journal of the J\'anos Bolyai
              Mathematical Society},
    VOLUME = {4},
      YEAR = {1984},
    NUMBER = {4},
     PAGES = {307--316},
      ISSN = {0209-9683},
   MRCLASS = {05C35},
  MRNUMBER = {779891},
       DOI = {10.1007/BF02579141},
       URL = {https://doi.org/10.1007/BF02579141},
    EPRINT = {},
EPRINTTYPE = {},
      NOTE = {},
}

@article {Thomason:1984:ExtremalContractions,
    AUTHOR = {Thomason, Andrew},
     TITLE = {An extremal function for contractions of graphs},
   JOURNAL = {Math. Proc. Cambridge Philos. Soc.},
  FJOURNAL = {Mathematical Proceedings of the Cambridge Philosophical
              Society},
    VOLUME = {95},
      YEAR = {1984},
    NUMBER = {2},
     PAGES = {261--265},
      ISSN = {0305-0041,1469-8064},
   MRCLASS = {05C35},
  MRNUMBER = {735367},
MRREVIEWER = {H.\ Kramer},
       DOI = {10.1017/S0305004100061521},
       URL = {https://doi.org/10.1017/S0305004100061521},
    EPRINT = {},
EPRINTTYPE = {},
      NOTE = {},
}

@article {Thomason:2001:CompleteMinors,
    AUTHOR = {Thomason, Andrew},
     TITLE = {The extremal function for complete minors},
   JOURNAL = {J. Combin. Theory Ser. B},
  FJOURNAL = {Journal of Combinatorial Theory. Series B},
    VOLUME = {81},
      YEAR = {2001},
    NUMBER = {2},
     PAGES = {318--338},
      ISSN = {0095-8956,1096-0902},
   MRCLASS = {05C83 (05C35)},
  MRNUMBER = {1814910},
MRREVIEWER = {W.\ Mader},
       DOI = {10.1006/jctb.2000.2013},
       URL = {https://doi.org/10.1006/jctb.2000.2013},
    EPRINT = {},
EPRINTTYPE = {},
      NOTE = {},
}

@article {AlonKS:2023:CompleteMinors,
    AUTHOR = {Alon, Noga and Krivelevich, Michael and Sudakov, Benny},
     TITLE = {Complete minors and average degree: a short proof},
   JOURNAL = {J. Graph Theory},
  FJOURNAL = {Journal of Graph Theory},
    VOLUME = {103},
      YEAR = {2023},
    NUMBER = {3},
     PAGES = {599--602},
      ISSN = {0364-9024,1097-0118},
   MRCLASS = {05C83 (05D40)},
  MRNUMBER = {4596514},
       DOI = {10.1002/jgt.22937},
       URL = {https://doi.org/10.1002/jgt.22937},
    EPRINT = {2202.08530},
EPRINTTYPE = {arxiv},
}

@unpublished{ChudnovskyGKM:2025:WheelSeparators,
    AUTHOR = {Chudnovsky, Maria and Gollin, J. Pascal and Krnc, Matja\v{z} and Milani\v{c}, Martin},
    TITLE = {Dominated balanced separators in wheel-induced-minor-free graphs},
    YEAR = {2025},
    NOTE = {Preprint},
    EPRINT = {2512.12329},
    EPRINTTYPE = {arxiv},
}

@book {Diestel:2025:GraphTheory,
    AUTHOR = {Diestel, Reinhard},
     TITLE = {Graph theory},
    SERIES = {Graduate Texts in Mathematics},
    VOLUME = {173},
   EDITION = {Sixth},
 PUBLISHER = {Springer, Berlin},
      YEAR = {2025},
     PAGES = {xx+454},
      ISBN = {978-3-662-70106-5},
   MRCLASS = {05-01 (05Cxx)},
  MRNUMBER = {4874150},
       URL = {https://diestel-graph-theory.com},
    EPRINT = {},
EPRINTTYPE = {},
}

@incollection {HickingbothamIMW:2024:Circlegraphs,
    AUTHOR = {Hickingbotham, Robert and Illingworth, Freddie and Mohar,
              Bojan and Wood, David R.},
     TITLE = {Treewidth, circle graphs and circular drawings},
 BOOKTITLE = {2021--2022 {MATRIX} annals},
    SERIES = {MATRIX Book Ser.},
    VOLUME = {5},
     PAGES = {213--227},
 PUBLISHER = {Springer, Cham},
      YEAR = {2024},
      ISBN = {978-3-031-47416-3},
   MRCLASS = {05C75 (05C62)},
  MRNUMBER = {4807277},
       DOI = {10.1007/978-3-031-47417-0\_12},
       URL = {https://doi.org/10.1007/978-3-031-47417-0_12},
    EPRINT = {2212.11436},
    EPRINTTYPE = {arxiv},
}

@thesis {Gartland:2023:Thesis,
       TYPE = {PhD Thesis},
     AUTHOR = {Gartland, Peter},
      TITLE = {Quasi-Polynomial Time Techniques for Independent Set and Beyond in Hereditary Graph Classes}, 
INSTITUTION = {UC Santa Barbara},
       YEAR = {2023},
        URL = {https://escholarship.org/uc/item/0kk6d2jv},
     EPRINT = {},
 EPRINTTYPE = {},
       NOTE = {},
}

@thesis{Aboulker:2013:Thesis,
       TYPE = {PhD Thesis},
     AUTHOR = {Aboulker, Pierre},
      TITLE = {Excluding slightly more than a cycle},
INSTITUTION = {Universit{\'e} Paris Diderot (Paris 7)},
       YEAR = {2013},
        URL = {https://www.di.ens.fr/~paboulker/content/these.pdf},
     EPRINT = {},
 EPRINTTYPE = {},
       NOTE = {},
}

@article {Gyarfas:1985:GyarfasPath,
    AUTHOR = {Gy\'arf\'as, Andr\'{a}s},
     TITLE = {Problems from the world surrounding perfect graphs},
 BOOKTITLE = {Proceedings of the {I}nternational {C}onference on
              {C}ombinatorial {A}nalysis and its {A}pplications
              ({P}okrzywna, 1985)},
   JOURNAL = {Zastos. Mat.},
  FJOURNAL = {Polska Akademia Nauk. Instytut Matematyczny. Zastosowania
              Matematyki},
    VOLUME = {19},
      YEAR = {1987},
    NUMBER = {3-4},
     PAGES = {413--441},
      ISSN = {0044-1899},
   MRCLASS = {05C15},
  MRNUMBER = {951359},
MRREVIEWER = {Alan\ Tucker},
    EPRINT = {https://users.renyi.hu/~gyarfas/Cikkek/36_Gyarfas_ProblemsFromTheWorldSurroundiPerfectGraphs.pdf},
EPRINTTYPE = {urltiny},
}

@article {ChudnovskyPPT:2024:ApproxMWIS,
    AUTHOR = {Chudnovsky, Maria and Pilipczuk, Marcin and Pilipczuk,
              Micha\l{} and Thomass\'e, St\'ephan},
     TITLE = {Quasi-polynomial time approximation schemes for the maximum
              weight independent set problem in {$H$}-free graphs},
   JOURNAL = {SIAM J. Comput.},
  FJOURNAL = {SIAM Journal on Computing},
    VOLUME = {53},
      YEAR = {2024},
    NUMBER = {1},
     PAGES = {47--86},
      ISSN = {0097-5397,1095-7111},
   MRCLASS = {68R10 (05C69 05C85 68W25)},
  MRNUMBER = {4708895},
MRREVIEWER = {Chuan-Min\ Lee},
       DOI = {10.1137/20M1333778},
       URL = {https://doi.org/10.1137/20M1333778},
    EPRINT = {1907.04585},
EPRINTTYPE = {arxiv},
      NOTE = {},
}

@article {DallardMS:2024:TWvsOmega1,
    AUTHOR = {Dallard, Cl\'ement and Milani\v c, Martin and \v Storgel,
              Kenny},
     TITLE = {Treewidth versus clique number. {I}. {G}raph classes with a
              forbidden structure},
   JOURNAL = {SIAM J. Discrete Math.},
  FJOURNAL = {SIAM Journal on Discrete Mathematics},
    VOLUME = {35},
      YEAR = {2021},
    NUMBER = {4},
     PAGES = {2618--2646},
      ISSN = {0895-4801,1095-7146},
   MRCLASS = {68R10 (05C05 05C40 05C69 05C75 05C83 05C85)},
  MRNUMBER = {4334541},
MRREVIEWER = {Mukkai\ S.\ Krishnamoorthy},
       DOI = {10.1137/20M1352119},
       URL = {https://doi.org/10.1137/20M1352119},
    EPRINT = {2006.06067},
EPRINTTYPE = {arxiv},
}

@article {DallardMS:2024:TWvsOmega2,
    AUTHOR = {Dallard, Cl\'ement and Milani\v c, Martin and \v Storgel,
              Kenny},
     TITLE = {Treewidth versus clique number. {II}. {T}ree-independence
              number},
   JOURNAL = {J. Combin. Theory Ser. B},
  FJOURNAL = {Journal of Combinatorial Theory. Series B},
    VOLUME = {164},
      YEAR = {2024},
     PAGES = {404--442},
      ISSN = {0095-8956,1096-0902},
   MRCLASS = {05C85 (05C69 05C75)},
  MRNUMBER = {4664382},
MRREVIEWER = {D.\ J.\ Guan},
       DOI = {10.1016/j.jctb.2023.10.006},
       URL = {https://doi.org/10.1016/j.jctb.2023.10.006},
    EPRINT = {2111.04543},
EPRINTTYPE = {arxiv},
}

@article {DallardMS:2024:TWvsOmega3,
    AUTHOR = {Dallard, Cl\'ement and Milani\v c, Martin and \v Storgel,
              Kenny},
     TITLE = {Treewidth versus clique number. {III}. {T}ree-independence
              number of graphs with a forbidden structure},
   JOURNAL = {J. Combin. Theory Ser. B},
  FJOURNAL = {Journal of Combinatorial Theory. Series B},
    VOLUME = {167},
      YEAR = {2024},
     PAGES = {338--391},
      ISSN = {0095-8956,1096-0902},
   MRCLASS = {05C75 (05C69 05C83)},
  MRNUMBER = {4730297},
MRREVIEWER = {D.\ J.\ Guan},
       DOI = {10.1016/j.jctb.2024.03.005},
       URL = {https://doi.org/10.1016/j.jctb.2024.03.005},
    EPRINT = {2206.15092},
EPRINTTYPE = {arxiv},
      NOTE = {},
}

@book {CyganFKLMPPS:2015:ParameterizedAlgorithms,
    AUTHOR = {Cygan, Marek and Fomin, Fedor V. and Kowalik, {\L}ukasz and
              Lokshtanov, Daniel and Marx, D\'aniel and Pilipczuk, Marcin
              and Pilipczuk, Micha\l{} and Saurabh, Saket},
     TITLE = {Parameterized algorithms},
 PUBLISHER = {Springer, Cham},
      YEAR = {2015},
     PAGES = {xviii+613},
      ISBN = {978-3-319-21275-3},
  MRNUMBER = {3380745},
       DOI = {10.1007/978-3-319-21275-3},
       URL = {https://doi.org/10.1007/978-3-319-21275-3},
    EPRINT = {},
EPRINTTYPE = {},
      NOTE = {},
}

@article{ChudnovskyGHLS:2026:EvenHoleFree,
    AUTHOR = {Chudnovsky, Maria and Gartland, Peter and Hajebi, Sepehr
             and Lokshtanov, Daniel and Spirkl, Sophie},
     TITLE = {Induced subgraphs and tree decompositions {XV}.
             {Even}-hole-free graphs with bounded clique number
             have logarithmic treewidth},
   JOURNAL = {Transactions of the American Mathematical Society},
      YEAR = {2026},
       DOI = {10.1090/tran/9807},
    EPRINT = {2402.14211},
EPRINTTYPE = {arxiv},
}

@article {CarmesinDHH:2014:kblocks,
    AUTHOR = {Carmesin, J. and Diestel, R. and Hamann, M. and Hundertmark,
              F.},
     TITLE = {{$k$}-blocks: a connectivity invariant for graphs},
   JOURNAL = {SIAM J. Discrete Math.},
  FJOURNAL = {SIAM Journal on Discrete Mathematics},
    VOLUME = {28},
      YEAR = {2014},
    NUMBER = {4},
     PAGES = {1876--1891},
      ISSN = {0895-4801,1095-7146},
   MRCLASS = {05C85 (05C05)},
  MRNUMBER = {3270977},
MRREVIEWER = {Litao\ Guo},
       DOI = {10.1137/130923646},
       URL = {https://doi.org/10.1137/130923646},
    EPRINT = {1305.4557},
EPRINTTYPE = {arxiv},
      NOTE = {},
}

@article {BellenbaumDiestel:2002,
    AUTHOR = {Bellenbaum, Patrick and Diestel, Reinhard},
     TITLE = {Two short proofs concerning tree-decompositions},
   JOURNAL = {Combin. Probab. Comput.},
  FJOURNAL = {Combinatorics, Probability and Computing},
    VOLUME = {11},
      YEAR = {2002},
    NUMBER = {6},
     PAGES = {541--547},
      ISSN = {0963-5483,1469-2163},
   MRCLASS = {05C05},
  MRNUMBER = {1940119},
       DOI = {10.1017/S0963548302005369},
       URL = {https://doi.org/10.1017/S0963548302005369},
    EPRINT = {},
EPRINTTYPE = {},
      NOTE = {},
}

@article {RobertsonST:1994:quicklyplanar,
    AUTHOR = {Robertson, Neil and Seymour, Paul and Thomas, Robin},
     TITLE = {Quickly excluding a planar graph},
   JOURNAL = {J. Combin. Theory Ser. B},
  FJOURNAL = {Journal of Combinatorial Theory. Series B},
    VOLUME = {62},
      YEAR = {1994},
    NUMBER = {2},
     PAGES = {323--348},
      ISSN = {0095-8956,1096-0902},
   MRCLASS = {05C10},
  MRNUMBER = {1305057},
MRREVIEWER = {Derek\ Smith},
       DOI = {10.1006/jctb.1994.1073},
       URL = {https://doi.org/10.1006/jctb.1994.1073},
    EPRINT = {},
EPRINTTYPE = {},
      NOTE = {},
}

@article {RobertsonSeymour:1986:GM5,
    AUTHOR = {Robertson, Neil and Seymour, P. D.},
     TITLE = {Graph minors. {V}. {E}xcluding a planar graph},
   JOURNAL = {J. Combin. Theory Ser. B},
  FJOURNAL = {Journal of Combinatorial Theory. Series B},
    VOLUME = {41},
      YEAR = {1986},
    NUMBER = {1},
     PAGES = {92--114},
      ISSN = {0095-8956,1096-0902},
   MRCLASS = {05C38 (68R10)},
  MRNUMBER = {854606},
       DOI = {10.1016/0095-8956(86)90030-4},
       URL = {https://doi.org/10.1016/0095-8956(86)90030-4},
    EPRINT = {},
EPRINTTYPE = {},
      NOTE = {},
}

@article {RobertsonSeymour:1995:GM13,
    AUTHOR = {Robertson, Neil and Seymour, P. D.},
     TITLE = {Graph minors. {XIII}. {T}he disjoint paths problem},
   JOURNAL = {J. Combin. Theory Ser. B},
  FJOURNAL = {Journal of Combinatorial Theory. Series B},
    VOLUME = {63},
      YEAR = {1995},
    NUMBER = {1},
     PAGES = {65--110},
      ISSN = {0095-8956,1096-0902},
   MRCLASS = {05C38 (05C85)},
  MRNUMBER = {1309358},
       DOI = {10.1006/jctb.1995.1006},
       URL = {https://doi.org/10.1006/jctb.1995.1006},
}

@article {AboulkerAKST:2021:evenholefreetw,
    AUTHOR = {Aboulker, Pierre and Adler, Isolde and Kim, Eun Jung and
              Sintiari, Ni Luh Dewi and Trotignon, Nicolas},
     TITLE = {On the tree-width of even-hole-free graphs},
   JOURNAL = {European J. Combin.},
  FJOURNAL = {European Journal of Combinatorics},
    VOLUME = {98},
      YEAR = {2021},
     PAGES = {Paper No. 103394, 21},
      ISSN = {0195-6698,1095-9971},
   MRCLASS = {05C75 (05C05)},
  MRNUMBER = {4286665},
MRREVIEWER = {Bo\v stjan\ Bre\v sar},
       DOI = {10.1016/j.ejc.2021.103394},
       URL = {https://doi.org/10.1016/j.ejc.2021.103394},
    EPRINT = {2008.05504},
EPRINTTYPE = {arxiv},
      NOTE = {},
}

@article {SintiariTrotignon:2021:LayeredWheels,
    AUTHOR = {Sintiari, Ni Luh Dewi and Trotignon, Nicolas},
     TITLE = {({T}heta, triangle)-free and (even hole, {$K_4$})-free
              graphs---part 1: {L}ayered wheels},
   JOURNAL = {J. Graph Theory},
  FJOURNAL = {Journal of Graph Theory},
    VOLUME = {97},
      YEAR = {2021},
    NUMBER = {4},
     PAGES = {475--509},
      ISSN = {0364-9024,1097-0118},
   MRCLASS = {05C75},
  MRNUMBER = {4313193},
MRREVIEWER = {Zden\v ek\ Ryj\'a\v cek},
       DOI = {10.1002/jgt.22666},
       URL = {https://doi.org/10.1002/jgt.22666},
    EPRINT = {1906.10998},
EPRINTTYPE = {arxiv},
      NOTE = {},
}

@unpublished{ChudnovskyTrotignon:2024:twvsomegacounterexamples,
    AUTHOR = {Chudnovsky, Maria and Trotignon, Nicolas},
    TITLE = {On treewidth and maximum cliques},
    YEAR = {2024},
    NOTE = {Preprint},
    EPRINT = {2405.07471},
    EPRINTTYPE = {arxiv},
}

@article {AbrishamiACHS:2024:ISTD7,
    AUTHOR = {Abrishami, Tara and Alecu, Bogdan and Chudnovsky, Maria and
              Hajebi, Sepehr and Spirkl, Sophie},
     TITLE = {Induced subgraphs and tree decompositions {VII}. {B}asic
              obstructions in {$H$}-free graphs},
   JOURNAL = {J. Combin. Theory Ser. B},
  FJOURNAL = {Journal of Combinatorial Theory. Series B},
    VOLUME = {164},
      YEAR = {2024},
     PAGES = {443--472},
      ISSN = {0095-8956,1096-0902},
   MRCLASS = {05C75},
  MRNUMBER = {4664383},
MRREVIEWER = {Ko-Wei\ Lih},
       DOI = {10.1016/j.jctb.2023.10.008},
       URL = {https://doi.org/10.1016/j.jctb.2023.10.008},
    EPRINT = {2212.02737},
EPRINTTYPE = {arxiv},
      NOTE = {},
}

@article {Bonnet:2025:SparseInducedSubgraphs,
    AUTHOR = {Bonnet, \'Edouard},
     TITLE = {Sparse induced subgraphs of large treewidth},
   JOURNAL = {J. Combin. Theory Ser. B},
  FJOURNAL = {Journal of Combinatorial Theory. Series B},
    VOLUME = {173},
      YEAR = {2025},
     PAGES = {184--203},
      ISSN = {0095-8956,1096-0902},
   MRCLASS = {05C75 (05C05)},
  MRNUMBER = {4881549},
MRREVIEWER = {Qiuzhenyu\ Tao},
       DOI = {10.1016/j.jctb.2025.03.002},
       URL = {https://doi.org/10.1016/j.jctb.2025.03.002},
    EPRINT = {2405.13797},
EPRINTTYPE = {arxiv},
      NOTE = {},
}

@article {ChuzhoyTan:2021:PolyGrid,
    AUTHOR = {Chuzhoy, Julia and Tan, Zihan},
     TITLE = {Towards tight(er) bounds for the excluded grid theorem},
   JOURNAL = {J. Combin. Theory Ser. B},
  FJOURNAL = {Journal of Combinatorial Theory. Series B},
    VOLUME = {146},
      YEAR = {2021},
     PAGES = {219--265},
      ISSN = {0095-8956,1096-0902},
   MRCLASS = {05C75 (05C83)},
  MRNUMBER = {4155282},
MRREVIEWER = {Benjamin\ Moore},
       DOI = {10.1016/j.jctb.2020.09.010},
       URL = {https://doi.org/10.1016/j.jctb.2020.09.010},
    EPRINT = {1901.07944},
EPRINTTYPE = {arxiv},
      NOTE = {},
}

\appendix

\section{}

In this short appendix we show, for completeness, that \zcref{conj:balancedseps} is equivalent to the conjecture of Gartland and Lokshtanov. The argument is essentially from~\cite{Gartland:2023:Thesis}. 

\begin{proposition}\label{prop:equivalent-dominated-separators}
The following statements are equivalent.
\begin{enumerate}
    \item \label{item:binary-weighted} For every planar graph~$H$, there exist a positive
    integer~$d$ and a~${c \in \halftoone}$ such that~$\Free{im}{H}$ has $d$-dominated binary-weighted $c$-balanced separators.

    \item \label{item:gartland-lokshtanov} For every planar graph~$H$, there exists a positive integer~$d$ such that~$\Free{im}{H}$ has $d$-dominated $\nicefrac{1}{2}$-balanced separators. 
\end{enumerate}
\end{proposition}

\begin{proof}
    Let us first prove that \ref{item:binary-weighted} implies \ref{item:gartland-lokshtanov}. 
    Take the constant weighting and repeatedly separate the unique component, if any, containing more than half the original vertices.
    Since the class is hereditary, after at most ${r = \lceil \log 2/\log(1/c)\rceil}$ applications we obtain an $rd$-dominated $\nicefrac12$-balanced separator. 

    So let us assume~\ref{item:gartland-lokshtanov} holds. 
    Fix a planar graph~$H$ and choose a grid~$Q$ containing~$H$ as an induced minor; such a grid exists since grids are universal graphs with respect to the induced minor relation; see, e.g.,~\cite[Theorem~12]{CampbellDDFGHHWWY:2024}. 
    We may assume that the minimum degree of~$Q$ is at least~$2$. 
    Let~$d$ be the constant supplied by the assumption for~$Q$.
    
    We use the pendant-vertex reduction described in~\cite[p.~46]{Gartland:2023:Thesis}.
    Let~$G$ be an $H$-induced-minor-free graph on~$n$ vertices, and let~${X \subseteq V(G)}$ be the set of vertices of weight~$1$. 
    Obtain~$G'$ by attaching~${n+1}$ pendant vertices to each vertex of~$X$.
    Then~$G'$ is $Q$-induced-minor-free: $G$ excludes~$Q$, and adding pendant vertices cannot create an induced minor of minimum degree at least two. 
    
    By assumption, $G'$ has a $\nicefrac{1}{2}$-balanced separator~$N_{G'}[S]$ with~${\abs{S} \leq d}$. 
    Replacing any pendant vertex in~$S$ by its neighbour only enlarges this closed neighbourhood, so we may assume~${S \subseteq V(G)}$.
    For every component~$C$ of~${G-N_G[S]}$, the vertices of~$C$ and all pendant vertices attached to them lie in one component of~$G'-N_{G'}[S]$. 
    Hence
    \[
        \abs{V(C)} + (n+1)\abs{V(C)\cap X}
        \leq \frac{n+(n+1)\abs{X}}{2},
    \]
    and therefore
    \[
        \abs{V(C) \cap X}
        \leq \frac{\abs{X}}{2} + \frac{n-2\abs{V(C)}}{2(n+1)}
        < \frac{\abs{X}+1}{2}.
    \]
    By integrality, we conclude that ${\abs{V(C)\cap X} \leq \abs{X}/2}$, proving \ref{item:binary-weighted}. 
\end{proof}

\end{document}